\documentclass[A4paper]{article}

\usepackage[T1]{fontenc}

\usepackage{amsthm,amssymb}
\usepackage{amsmath}
\usepackage{mathtools}
\usepackage{bm}
\usepackage{amsfonts}
\usepackage{subfiles}
\usepackage{multicol}
\usepackage{tabularx}
\usepackage{subcaption}
\usepackage{caption}
\usepackage[absolute,overlay]{textpos}
\usepackage{ifthen}
\usepackage{mathrsfs}
\usepackage{scalerel}
\usepackage{authblk}

\usepackage[]{graphicx}
\usepackage{tikz}
\usetikzlibrary{arrows,positioning}
\usetikzlibrary{decorations.markings}
\tikzset{myptr/.style args={#1}{%
    decoration={%
      markings,%
      mark=at position #1 with {\arrow[scale=2,>=latex]{>}}%
    },%
    postaction={decorate},%
    shorten >=8pt%
  }%
}
\usetikzlibrary{patterns}
\usetikzlibrary{decorations.pathmorphing}
\usetikzlibrary{calc}
\colorlet{tikzcolorA}{red!80!black}
\colorlet{tikzcolorB}{black}
\colorlet{tikzcolorC}{green!70!black}
\colorlet{tikzcolorD}{gray!60!white}

\usepackage[most]{tcolorbox}
\tcbuselibrary{listings,breakable}

\definecolor{moplorange}{rgb}{0.96, 0.82, 0.66}
\definecolor{moporange}{rgb}{0.92, 0.6, 0.28}
\definecolor{moplblue}{rgb}{0.10, 0.48, 0.86}
\definecolor{mopblue}{rgb}{0.07, 0.35, 0.63}
\definecolor{mopdblue}{rgb}{0.05, 0.25, 0.45}
\definecolor{mygreen}{rgb}{0.058,.596,0.047}  

\usepackage{hyperref}
\hypersetup{
     colorlinks   = true,
     linkcolor    = blue,
     citecolor    = red
}

\tcbuselibrary{theorems}
\tcbuselibrary{skins}

\usepackage{algorithm}
\usepackage{algpseudocode}

\usepackage{mycommands}
\usepackage{xcolor}

\newcommand{\opspd}[1]{\mathcal S_{#1}}  
\newcommand{\mprox}[2][]{\mathrm{prox}_{#2}^{#1}}

\DeclareMathOperator{\iconv}{%
\mathbin{\ooalign{\raisebox{-.5ex}{$\triangledown$}\cr\hidewidth\rule[.3ex]{.4pt}{1ex}\hidewidth}}}%
\newtheorem{theorem}{Theorem}[section]

\newtheorem{lemma}[theorem]{Lemma}

\newtheorem{assum}{Assumption}
\newtheorem{remark}{Remark}

\theoremstyle{definition}

\usepackage[colorinlistoftodos,prependcaption,textsize=tiny]{todonotes}

\newcommand{\appendixhead}%
{\textbf{Supplement Material: A Quasi-Newton Primal-Dual Algorithm with Line Search}
\vspace{0.25in}}
\begin{document}
\title{A Quasi-Newton Primal-Dual Algorithm with Line Search\thanks{ We acknowledge funding by the ANR-DFG joint project TRINOM-DS under the number DFG OC150/5-1.
  }}
%
%
\author[*]{Shida Wang}
\author[**]{Jalal Fadili}
\author[*]{Peter Ochs}
\affil[*]{Department of Mathematics, University of Tübingen, Germany}
\affil[**]{Normandie Universit\'e, ENSICAEN, UNICAEN, CNRS, GREYC, France.}

%
%
\maketitle              
\begin{abstract}
    Quasi-Newton methods refer to a class of algorithms at the interface between first and second order methods. They aim to progress as substantially as second order methods per iteration, while maintaining the computational complexity of first order methods. The approximation of second order information by first order derivatives can be expressed as adopting a variable metric, which for (limited memory) quasi-Newton methods is of type ``identity $\pm$ low rank''. This paper continues the effort to make these powerful methods available for non-smooth systems occurring, for example, in large scale Machine Learning applications by exploiting this special structure. We develop a line search variant of a recently introduced quasi-Newton primal-dual algorithm, which adds significant flexibility, admits larger steps per iteration, and circumvents the complicated precalculation of a certain operator norm. We prove convergence, including convergence rates, for our proposed method and outperform related algorithms in a large scale image deblurring application. 
\end{abstract}

\section{Introduction}

In modern optimization, the datasets and dimensionality of the problems and parameters is vastly increasing. In the early stages of large scale optimization, the shift from second order to first order optimization could cope with the increasing dimensionality of the problems. Second order methods achieve a significant progress per iteration at the cost of a high computational load, since the computation of the second derivative (Hessian) and oftentimes its inverse are required, which is intractable in the large scale regime. In contrast, first order methods have a low computational effort per iteration but also less information about the objective. The gradient cannot capture curvature information and hence may fail to provide directions that allow for large steps. Nevertheless, for the large scale regime this exchange pays off. 

However, the ever increasing dimensionality of the considered problems and datasets asks for faster algorithms. Motivated by classical optimization and the discussion above, algorithms at the interface of first and second order methods are key to reach the next level. Tractability in the (nowadays extremely) large scale regime requires methods that are as cheap as first order methods, while progressing as substantially as second order algorithms: \emph{Quasi-Newton methods}. While they are known for their outstanding performance in unconstrained smooth optimization, their development for non-smooth (or constrained smooth) problems is insufficiently understood. As we discuss in related work below, most algorithmic development is either too simplistic, in the sense that only a diagonal metric is admitted, which can hardly capture second order information of the objective, or too theoretical, in the sense that a good performance is proved in theory while the implementation cost is on a par with that of second order methods. Algorithmic subproblems (e.g. the evaluation of the proximal mapping) that are easy to solve with respect to the Euclidean metric may become intractable with respect to another metric.

We pursue the line of research initiated in \cite{becker2012quasi,becker2019quasi} that considers both aspects as equally important. Key is the observation that quasi-Newton methods actually generate a metric of the specific type ``identity $\pm$ low rank'', which allows for an efficient proximal calculus (cf. Section~\ref{sec:prox-calc}) that unlocks the quasi-Newton power---well-known from classical optimization---in the area of optimization for Machine Learning. This idea was recently transferred to non-convex optimization in \cite{kanzow2021globalized,kanzow2022efficient} and to monotone inclusion problems in \cite{WFO22}. A special case of the latter setting comprises the extremely broad class of convex--concave saddle point problems, which has numerous applications in Machine Learning, Computer Vision, Image Processing and Statistics \cite{chambolle2016ergodic,chambolle2016introduction,goldstein2015adaptive,valkonen2014primal,applegate2021practical}. 

In this paper, we restrict our interest to saddle-point problems only. This focus allows us to design a quasi-Newton primal--dual algorithm that is tailored to this setting and therefore highly efficient and adaptable thanks to an additional line search procedure. This has several advantages as compared to a fixed step size: (i) the oftentimes expensive computation of the operator norm can be avoided, (ii) the choice of metric need not obey any static spectral restrictions, and (iii) in many situations larger steps and thus a faster convergence is observed. 

Our main contribution is reduction of the gap between the outstanding performance of quasi-Newton methods in classical optimization and quasi-Newton methods for (non-smooth) convex--concave saddle point problems for modern optimization in Machine Learning. In detail, our contribution is the following:
\begin{enumerate}
    \item We extend the line-search based primal--dual algorithm in \cite{malitsky2018first} (extension of \cite{chambolle2011first} by line search) to incorporate a quasi-Newton metric with efficiently implementable proximal mapping; thereby aiming equally at theoretical convergence guarantees (including convergence rates) as well as highly efficient implementation.
    \item We unlock the use of multi-memory quasi-Newton metrics (L-BFGS and SR1 method) via a compact representation for primal--dual algorithms, including their efficient implementation via a semi-smooth Newton solver.
    \item The proposed algorithm outperforms the line-search based primal--dual algorithm (with identity metric) on a challenging image deblurring problem.
\end{enumerate}

\subsection{Related Work}

Due to page limitations, for an extended discussion of quasi-Newton approaches in non-smooth optimization and the vast literature on primal-dual algorithms, we refer to \cite{WFO22}.

\textbf{Non-smooth quasi-Newton.} For a class of non-smooth problems that are given as a composition of a smooth function $h$ and a non-smooth function $g$, \cite{patrinos2014forward,stella2017forward}
combine quasi-Newton methods with forward--backward splitting via the forward--backward envelope. If $g$ is an indicator function, \cite{schmidt,schmidt201211} proposed a projected quasi-Newton method which requires either solving a complicated subproblem or is restricted to a diagonal metric. Later, their work was extended by \cite{lee2014proximal} to a more general setting. \cite{Vavasis,becker2012quasi,becker2019quasi} developed algorithms with efficient evaluation of the proximal operator with respect to a low-rank perturbed metric. It is worth to mention that in \cite{becker2012quasi,becker2019quasi} the subproblem is a low dimensional root finding problem which can be solved efficiently. Inspired by \cite{becker2019quasi}, the work in \cite{kanzow2022efficient} applied the limited-memory quasi-Newton method on non-convex problems. 

\textbf{PDHG.} Primal-Dual Hybrid Gradient (PDHG) is widely used to solve saddle point problems \cite{chambolle2011first,chambolle2016ergodic}. However, in order to guarantee the convergence of PDHG, the computation of the norm of a operator $K$ is required. To avoid this disadvantage, \cite{malitsky2018first} combined line search with PDHG to get a new algorithm PDAL. For faster convergence, variable metric is being used \cite{goldstein2015adaptive}. It shows the potential of combining quasi-Newton methods and PDAL via a variable metric, however suffers again from the need to solve more complicated subproblems, which is remedied in \cite{WFO22} for the more general class of monotone inclusion problems and hence builds the grounding for our proposed line search variant.


\section{Problem Setup: A class of saddle point problems}

Let $X$ and $Y$ be finite dimensional real vector spaces with inner product $\scal{\cdot}{\cdot}$ and induced norm $\vnorm\cdot=\sqrt{\scal\cdot\cdot}$. 
We consider saddle point problems
\begin{equation}\label{eq:saddle-point-problem}
    \min_{x\in X} \max_{y\in Y}\, \scal{Kx}{y} + g(x) + h(x) -f^*(y)\,,
\end{equation}
where $\map{g,h}{X}{\eR}$ are proper, lower semi-continuous (lsc), convex functions with $h$, additionally, having an $L$-Lipschitz continuous gradient, $\map{f^*}{Y}{\eR}$ is a proper, lsc, convex function; the convex conjugate (Legendre--Fenchel conjugate) of a function $f$, and $\map{K}{X}{Y}$ is a bounded linear operator with operator norm $L_K:=\norm{K}$ and adjoint $K^*$. Moreover, throughout this paper we assume that \eqref{eq:saddle-point-problem} has a saddle point. 
We remark that \eqref{eq:saddle-point-problem} is equivalent to the \emph{primal problem}
\begin{equation}\label{eq:primal-problem}
    \min_{x\in X}\, f(Kx) + g(x) + h(x) \,,
\end{equation}
and to the \emph{dual problem}
\begin{equation}\label{eq:dual-problem}
    \max_{y\in Y}\, -\Big(f^*(y) + (g^*\iconv h^*)(-K^* y)\Big)\,,
\end{equation}
where we use the fact that the conjugate of a sum of two function $(g+h)^*$ equals the infimal convolution of the conjugate functions $g^*\iconv h^*$.

\section{Our quasi-Newton Primal-Dual Algorithm with Line Search}

The primal--dual algorithm that we develop in this paper is an extension of \cite{malitsky2018first} to incorporate a variable metric of quasi-Newton type, which itself adds an efficient line search procedure to the primal--dual hybrid gradient (PDHG) algorithm \cite{chambolle2011first} (aka Chambolle--Pock algorithm) and the extension in \cite{lorenz2015inertial}. The handling of (a possibly) non-smooth functions essentially relies on evaluating the proximal mapping, which we define here for a function $g$ and parameter $\tau$ with respect to an arbitrary metric $M\in \opspd{\alpha}(X)$, $\alpha>0$. Here, $\opspd{\alpha}$ is the set of bounded self-adjoint linear operators from Hilbert space $X$ to $X$ such that $M-\alpha\opid$ is positive semi-definite for each $M\in\opspd{\alpha}$. For simplicity, some notations are introduced:
\[
\norm[M]{x}^2 = \scal{Mx}{x},\quad x\in X\,.
\]
\[
   \mprox[M]{\tau g}(\bar x) := \argmin_{x\in X}\, g(x) + \frac{1}{2\tau}\vnorm[M]{x-\bar x}^2\,,\quad \text{and set}\ \prox{\tau}{g} := \mprox[\opid]{\tau g} \,,
\]
where $\opid$ is the identity (Euclidean) metric.

Algorithm~\ref{Alg:QN-PDHG-line} presents the proposed algorithm.
\begin{algorithm}[th]
\caption{Quasi-Newton PDHG with Line Search}\label{Alg:QN-PDHG-line}
\begin{algorithmic}
\Require [initial data] $x^1\in X$, [initial data] $y^0\in Y$, [maximal iteration count] $N \geq 0$, [scaling of line search parameter] $\mu\in(0,1)$, [extrapolation parameter] $\theta_0$, [initial dual step size] $\sigma_{0}$, [tolerance weight] $\delta\in(0,1)$, [primal-dual step ratio] $+\infty>\beta\geq\beta_{k}\geq\beta_{k+1}>0$, $\forall \k\in\N$.
\State \textbf{Update for $k=1,\ldots,N$}:
\begin{enumerate}
    \item[(i)] Compute $M_k$ according to a quasi-Newton framework (cf. Section~\ref{sec:bfgs-sr1-metric}).
    \item[(ii)] Compute dual update step:
    \begin{equation}\label{eq:dual-update}
        y^{\k} =\mathrm{prox}^{}_{\sigma_{k-1} f^*}( y^{k-1} +\sigma_{k-1}Kx^k)\,.
    \end{equation}
    \item[(iii)] Select
    $
        \bar\sigma_k \in [\frac{\beta_{k-1}}{\beta_k}\sigma_{k-1}, \sqrt{(1+\theta_{k-1})}\frac{\beta_{k-1}}{\beta_k}\sigma_{k-1}]
    $ 
    
    and compute the quasi-Newton primal update step by:
        \item[] \textbf{Line search:} Find the smallest power $i=0,1,2,\ldots$ such that 
              \begin{equation}\label{eq:primal-update}
              \begin{split}
                  \bar y^{k} =&\ y^{k} + \theta_k(y^{k}-y^{k-1})\,,\\
                  x^{k+1} =&\ \mathrm{prox}^{M_k}_{\tau_k g}\big( x^k-\tau_kM_k^{-1}K^* \bar y^{k}-\tau_k M_k^{-1}\nabla h( x^k)\big)
              \end{split}
              \end{equation}
              with  
              \[ 
               \sigma\piter\k = \bar\sigma\piter\k\cdot \mu^i\,, \quad 
                  \theta_k = \frac{\sigma_k}{\sigma_{k-1}}\,, \quad\text{and}\quad 
                  \tau_k = \beta_k\sigma_k
              \]
              satisfy
              \begin{multline}\label{CB}
                      \tau_k\sigma_k\norm[]{Kx^{k+1}-Kx^k}^2+ 2\tau_k \Big(h(x^{k+1})-h(x^k)-\scal{\nabla h(x^k)}{x^{k+1}-x^k}\Big)\\
                      \leq \delta \norm[M_k]{x^{k+1}-x^k}^2\,.
              \end{multline}
\end{enumerate}
\State \textbf{End of for-loop}
\end{algorithmic}
\end{algorithm}
The algorithm alternates between updates of the dual variable \eqref{eq:dual-update} and the primal variable \eqref{eq:primal-update}, where the line search is only implemented in the primal update. Let us discuss the algorithm for a fixed iteration $k$, i.e., we are given $x^k$, $y\iter\km$, $\sigma_{k-1}$, $\theta_{k-1}$ and a monotone decreasing sequence of $\beta_k$. The discussion of the quasi-Newton type variable metric in Step~(i) is deferred to Section~\ref{sec:bfgs-sr1-metric}. Step~(ii) is a standard dual update step. In Step~(iii), we perform the line search. We select a basic step size $\bar \sigma\piter\k$ and, in the $i$th loop of the line search, we perform a trial step \eqref{eq:primal-update} with the current $\sigma\piter\k=\bar\sigma\piter\k\cdot\mu^i$ and check if the breaking condition \eqref{CB} is satisfied. If `yes', the current iteration $k$ is completed. If `no', the new trial step size is reduced to $\sigma\piter\k=\bar\sigma\piter\k\cdot\mu^{i+1}$ (in the subsequent $(i+1)$th line search step).
Here, $M_k\in\mathcal{S}_\alpha(X)$ is symmetric positive definite and $\alpha \in (0,1)$.If $M_k$ is chosen as an identity, then we recover the breaking (stopping) criterion used in \cite{malitsky2018first}. However, in this paper, we adopt a variable metric $M_k$ which is generated by quasi-Newton methods to exploit the local geometry of the function $h$. As a result, it is more likely to obtain a larger step size $\sigma_k$ and fewer inner loops for the line search procedure as compared to the Euclidean version ($M_k=\opid$). The employed metric is of type ``identity $\pm$ low rank'' for which the proximal mapping can be computed efficiently as shown in Section~\ref{sec:prox-calc}.
%
\begin{remark}
    While the line search procedure is formulated for the primal problem, by duality, the primal problem can be interpreted as the dual of the dual problem and, thus, the dual problem as the primal problem. As a consequence, an equivalent algorithm with line search on the dual can be easily stated.
\end{remark}

\textbf{Discussion of computational cost for line search.} In general, every loop of the line search procedure requires recomputing several quantities, including \eqref{eq:primal-update} and $Kx\iter\kp$, $h(x\iter\kp)$ and $\norm[M_k]{x\iter\kp-x\iter\k}^2$ in \eqref{CB}. While this seems to be expensive at first glance, often \eqref{CB} is satisfied after 1--3 trial steps and hence large steps are taken with a low cost, as we underline in our experiments in Section~\ref{sec:experiment}. Nevertheless, the cost can be further reduced significantly in certain special cases, observed in \cite{malitsky2018first} and generalized here to our setting, whenever $\prox{\tau_k}{g}^{M_k}$ is a linear or affine operator.
\begin{enumerate}
    \item If $g(x) = \scal{c}{x}$, then $\prox{\tau_k}{g}^{M_k}(u) = u-\tau_k M_k^{-1}c$ and therefore, we obtain
    \begin{equation*}
    \begin{split}
         x^{k+1}&=\prox{\tau_k}{g}^{M_k}(x^{k}-\tau_k M_k^{-1}K^* \bar{y}^{k}-\tau_kM_k^{-1}\nabla h(x^{k})) \\&= x^{k}-\tau_k[M_k^{-1}K^* \bar{y}^{k}+M_k^{-1}\nabla h(x^{k})   + M_k^{-1}c]\,.
    \end{split}
    \end{equation*}
    \begin{equation*}
        Kx^{k+1}=Kx^{k}-\tau_k[KM_k^{-1}K^* \bar{y}^{k}+\tau_kKM_k^{-1}\nabla h(x^{k})+\tau_k KM_k^{-1}c]\,. 
    \end{equation*}
    \item If $g(x) = \frac{1}{2}\norm[]{x-b}^2$, then $\prox{\tau_k}{g}^{M_k}(u) = (\opid+\tau_kM_k^{-1})^{-1}[u+\tau_kM_k^{-1}b]$
    and therefore, we obtain
    \begin{equation*}
    \begin{split}
        x^{k+1} &= \prox{\tau_k}{g}^{M_k}(x^{k}-\tau_kM_k^{-1}K^* \bar{y}^{k}-\tau_kM_k^{-1}\nabla h(x^{k})) \\& = (\opid+\tau_kM_k^{-1})^{-1}[ x^{k}-\tau_kM_k^{-1}K^* \bar{y}^{k}-\tau_kM_k^{-1}\nabla h(x^{k}) +\tau_kM_k^{-1}b]\,,
    \end{split}
    \end{equation*}
    \begin{equation*}
        Kx^{k+1}=K(\opid+\tau_kM_k^{-1})^{-1}[x^k- \tau_k (M_k^{-1}K^* \bar{y}^{k}+M_k^{-1}\nabla h(x^{k})-M_k^{-1}b)] \,. 
    \end{equation*}
    \item Let $g(x) = \delta_H(x)$, where $H$ refers to the hyperplane $H:= \{u:\scal{u}{a}=b\}$. Then 
    $\prox{\tau_k}{g}^{M_k}(u)= u+ \frac{b-\scal{u}{a}}{\norm[M_k^{-1}]{a}^2}M_k^{-1}a$ and therefore, we obtain
    \begin{equation*}
        \begin{split}
            x^{k+1}&=\prox{\tau_k}{g}^{M_k}(x^{k}-\tau_kM_k^{-1}K^* \bar{y}^{k}-\tau_kM_k^{-1}\nabla h(x^{k}))
            \\&= x^{k}-\tau_k[M_k^{-1}K^* \bar{y}^{k}+M_k^{-1}\nabla h(x^{k})]\\&+ \frac{b-\scal{x^{k}-\tau_k[M_k^{-1}K^* \bar{y}^{k}+M_k^{-1}\nabla h(x^{k})]}{a}}{\norm[M_k^{-1}]{a}^2}M_k^{-1}a\,,\\
        \end{split}
    \end{equation*}
    \begin{equation*}
    \begin{split}
        Kx^{k+1} &= Kx^k -\tau_k[KM_k^{-1}K^* \bar{y}^{k}+KM_k^{-1}\nabla h(x^{k})] \\&+ \frac{b-\scal{x^k- \tau_k[M_k^{-1}K^* \bar{y}^{k}+\tau_kM_k^{-1}\nabla h(x^{k})]}{a}}{\norm[M_k^{-1}]{a}^2}KM_k^{-1}a \,.
    \end{split}
    \end{equation*}
\end{enumerate}

\section{Convergence Analysis of Algorithm~\ref{Alg:QN-PDHG-line}}

Let us now analyze the convergence of Algorithm~\ref{Alg:QN-PDHG-line}. As for most variable metric primal--dual algorithms (cf. \cite{davis2015convergence,combettes2014variable,combettes2014forwardbackward,WFO22}), we require the following restriction for the change of the metric from one iteration to the next. Under this condition, we can generalize all convergence results from \cite{malitsky2018first} by adapting their proofs.
\begin{assum}\label{assumption1} Let $\alpha\in (0,+\infty)$.
$(M_k)_{k\in\mathbb{N}}$ is a sequence in $\opspd{\alpha}(X)$ such that
\begin{equation}
\begin{cases}
    \exists C_M\in \R, \mathrm{s.t.}\mathrm{sup}_{k\in\mathbb{N}}||M_k||\leq C_M<\infty\,,\\
    (\exists(\eta_k)_{k\in\mathbb{N}}\in \ell_+^1(\mathbb{N}))(\forall k\in\mathbb{N})\colon\quad(1+\eta_k)M_{k}\succeq M_{k+1}\,.\\
    \end{cases}
\end{equation}
\end{assum}
\begin{lemma}\label{lemma:sigma}
    \begin{enumerate} 
        \item[(i)] There exists some $\sigma>0$ such that $\sigma_k\geq \sigma$ for any $k\in\N$.
        \item[(ii)] The line search terminates.
        \item[(iii)] If $\beta_k\equiv \beta$, $\theta_k$ is bounded from above by some $\theta$ for any $k\in\N$.
    \end{enumerate}
\end{lemma}
The proof is provided in Section~\ref{sec:proof-lemma}.
\begin{theorem}\label{thm:main-conv}
Consider Problem \eqref{eq:saddle-point-problem} and let the sequence $(x^k,y^k)_{k\in\N}$ be generated by Algorithm~\ref{Alg:QN-PDHG-line} with $\beta_k\equiv\beta$ where Assumption~\ref{assumption1} holds. Then $(x^k,y^k)_{k\in\N}$ is a bounded sequence and its cluster points are solutions of \eqref{eq:saddle-point-problem}. Furthermore, if $f^*\vert_{\mathrm{dom} f^*}$ is continuous and $\sigma_k$ is bounded from above, then the whole sequence $(x^k,y^k)_{k\in\N}$ converges to a solution of \eqref{eq:saddle-point-problem}.
\end{theorem}
The proof is provided in Section~\ref{sec:proof-main-conv}. \\

We obtain the same ergodic convergence rate as in \cite{malitsky2018first}, with respect to the primal--dual gap $\mathcal{G}_{\hat x,\hat y}$ which is the difference (gap) between the optimal primal objective value in \eqref{eq:primal-problem} and the optimal dual objective value \eqref{eq:dual-problem}.
\begin{theorem}\label{thm:conv-rate}
Let the sequence $(x^k,y^k)_{k\in\N}$ be generated by Algorithm \ref{Alg:QN-PDHG-line} with $\beta_k\equiv\beta$ where Assumption \ref{assumption1} holds and $(\hat x,\hat y)$ be some saddle point of \eqref{eq:saddle-point-problem}. Then it holds for a constant $D =\Pi_{k\in\N}(1+\eta_k)<+\infty$ that
\begin{equation}\label{eq:thm:conv-rate}
    \mathcal{G}_{\hat x,\hat y}(\bar X^N,\bar Y^N)\leq \frac{D}{s_N}\Big( \frac{1}{2\beta}\norm[M_1]{x^1-\hat x}^2 + \frac{1}{2}\norm[]{y^1-\hat y}^2+\sigma_1\theta_1D_{\hat x,\hat y}(y^0)\Big) = O\Big(\frac{1}{N}\Big)\,,
\end{equation}
where $s_N := \sum_{k=1}^N \sigma_k$, $ \bar X^N :=\frac{\sum_{k=1}^N \sigma_k x^{k+1}}{s_N} $ and $ \bar Y^N := \frac{\sigma_1\theta_1y^0+\sum_{k=1}^N\sigma_k \bar y^k}{\sigma_1\theta_1 +s_N}$. The last equality in \eqref{eq:thm:conv-rate} provides a simplified rate in terms of the big-O notation.
\end{theorem}
The proof is provided in Section~\ref{sec:proof-conv-rate}.

Under the additional assumption that $g$ is strongly convex, improved convergence rates can be derived when $\seq[\k\in\N]{\beta\piter\k}$ is varied appropriately. 
\begin{theorem}\label{thm:conv-rate-strong-convex}
  Assume $g$ is $\gamma$-strongly convex and $(x^k,y^k)_{k\in\N}$ is generated by Algorithm~\ref{Alg:QN-PDHG-line} with 
  \begin{equation}
      \beta_k=\frac{\beta_{k-1}}{\min\{1+\frac{\gamma}{C_M}\beta_{k-1}\sigma_{k-1},C_\theta\}}\,,\ \forall \k\in\N\,, \quad \text{and}\quad \beta_0 >0 \,,
  \end{equation}
  where $C_\theta\in \R_+$ is a constant, Assumption \ref{assumption1} holds and $(\hat x, \hat y)$ be some saddle point of \eqref{eq:saddle-point-problem}. Then, we have $(\theta_k)_{k\in\N}$ is bounded from above. Furthermore, we obtain
  \[
    \norm[]{x^N-\hat x}=O(1/N)\quad\text{and}\quad \mathcal{G}_{\hat x,\hat y}(\bar X^N,\bar Y^N)=O(1/N^2)\,,
  \]
  where $(\bar X^N,\bar Y^N)$ are the ergodic sequences defined in Theorem~\ref{thm:conv-rate}.
\end{theorem}
The proof is provided in Section~\ref{sec:proof-conv-rate-strong-convex}.
\begin{remark}
    For the result in Theorem~\ref{thm:conv-rate-strong-convex}, $\delta=1$ is also admitted.
\end{remark}
\section{Computing and Representing the quasi-Newton Metric} \label{sec:bfgs-sr1-metric}

In this section, we abuse notation in order to follow the conventions of quasi-Newton methods\footnote{For example, the variable $y\iter\k$ defined in \eqref{eq:s-and-y} is not the dual variable in Algorithm~\ref{Alg:QN-PDHG-line}.}. The metric $M_k$ is expected to be an approximation of the Hessian $\nabla^2 h(x^k)$ at $x^k$ for the $k$-th iteration. The most popular quasi-Newton methods are the SR1 and BFGS methods (and their low-memory variants), which update $M_k$ by adding a rank-one modification (SR1 method)
\begin{equation}
    M_{k+1}:= M^{SR1}_{k+1}=M_k + \frac{ (y^k-M_ks^k)(y^k-M_ks^k)^\top }{ (y^k-M_ks^k)^\top s^k}\,,
\end{equation}
or a rank-two modification (BFGS method)
\begin{equation}
    M_{k+1}\coloneqq M^{BFGS}_{k+1}= M_k + \frac{y^k(y^k)^\top }{(s^k)^\top y^k } - \frac{M_ks^k(s^k)^\top M_k}{(s^k)^\top M_ks_k}\,,
\end{equation}
respectively, where
\begin{equation}\label{eq:s-and-y}
    s^k:=x^{k+1}-x^k\quad\text{and}\quad y^k := \nabla h(x^{k+1}) - \nabla h(x^k) \,.
\end{equation}
In order to apply quasi-Newton methods on large-scale problems, \emph{$m$-limited memory quasi-Newton methods} are adopted \cite{kanzow2022efficient}, with the most popular version being L-BFGS \cite{wright1999numerical}. It means that instead of generating $M_k$ via all previous $s^i$ and $y^i$ for $i=1,\ldots, k$ and $M_0$, for each $k$, the metric $M_k$ is re-computed based on $M_{k,0}$ and the most recent $m$ vectors $s^i$ and $y^i$ for $i=k-m+1,\ldots,k$, if $k\geq m$. Usually, $m$ is very small, such that only a small storage will be required. 
As pointed out by \cite{byrd1994representations}, the matrices of the $m$-limited memory version of quasi-Newton methods have a compact representation of the form
\begin{equation}\label{eq:metric-compact}
    M_k=M_{k,0}+ A_kQ_k^{-1}A_k^\top \,,
\end{equation}
where $M_{k,0}\in\R^{n\times n}$, $n=\dim(X)$, is a symmetric positive definite matrix, $A_k\in \R^{n\times m}$, and a symmetric and non-singular matrix $Q_k\in\R^{m\times m}$ ($m\ll n$). 
For limited memory BFGS (known as L-BFGS), we have the following block-matrix representation
\begin{equation}
    A_k = \begin{bmatrix}M_{k,0}S_k&Y_k\end{bmatrix}\in \R^{n\times 2m}\quad\mathrm{and}\quad Q_k=\begin{bmatrix}
        -S_k^* M_{k,0}S_k&-L_k\\
        -L_k^* & D_k\\ 
    \end{bmatrix}\in \R^{2m\times 2m}\,,
\end{equation}
where $S_k$ and $Y_k$ are matrices collecting the $m$ most recent vectors in \eqref{eq:s-and-y} as columns, $D_k := D(S_k^\top Y_k)$ and $L_k :=L(S_k^\top Y_k)$ refer to the diagonal $D(\cdot)$ and the strict lower triangular $L(\cdot)$ part of the matrix $S_k^\top Y_k$, respectively. By using a spectral decomposition $Q^{-1}=V\Lambda V^\top $ with orthogonal $V\in \R^{s\times s}$ and diagonal $\Lambda\in\R^{s\times s}$, for some $s\in\N$, \eqref{eq:metric-compact} is transformed into the compact representation
\begin{equation}\label{metric}
    M_k = M_{k,0} + U_1U_1^\top  - U_2 U_2^\top \,,
\end{equation}
for some $U_1\in\R^{n\times m}$ and $U_2\in\R^{n\times m}$. In detail, since $\Lambda$ is a diagonal matrix with eigenvalues $\lambda_i$, $i=1,2,\ldots,s$ of $Q_k^{-1}$ on the diagonal, we decompose $\Lambda = \Lambda_1 -\Lambda_2$ where $\Lambda_1$, given by $(\Lambda_1)_{i,i}=\max(\lambda_i,0)$, $i=1,2,\ldots,s$, corresponds to the positive eigenvalues  and $\Lambda_2$, given by $(\Lambda_2)_{i,i}=\max(-\lambda_i,0)$, $i=1,2,\ldots,s$, corresponds to the negative eigenvalues. In this way, we obtain 
\begin{equation}
    U_1\coloneqq (A_kV)\Lambda_1^{1/2}\quad\mathrm{and}\quad U_2\coloneqq (A_kV)\Lambda_2^{1/2}\,.
\end{equation}
Theoretically, it is guaranteed that $M_k=M_{k,0}+ U_1U_1^\top-U_2U_2^\top $ is positive definite \cite{fletcher2013practical} if $s^ky^k>0$ for any $k\in\N$. However, in order to account for numerical rounding errors and the assumption that $M_k\in\opspd{\alpha}$ is bounded from above by some $C_M$, we adopt a scaling version:

\begin{equation}\label{eq:metric-scalling}
\begin{split}
    \tilde M_k & = M_{k,0}+ \gamma_1U_1U_1^\top - \gamma_2U_2U_2^\top \,,\\
    M_k &= \min\{\frac{C_M-\alpha}{\norm[2]{\tilde{M}_k}} ,1\}\tilde M_k +\alpha\opid \,,\\
\end{split}
\end{equation}

where $\norm[2]{\tilde M_k}$ denotes the $l_2$ norm of matrix $\tilde M_k$ and we set $\alpha=0.01$,$\gamma_1= 1$, $\gamma_2 = 1, C_M = 50$ in practice. There is an easy way to make sure that Assumption~\ref{assumption1} is satisfied by setting $\gamma_1 = \frac{\eta_k}{\norm[]{U_1}^2}$ and  $\gamma_2 = \frac{\eta_k}{\norm[]{U_2}^2}$ with arbitrary $\eta_k\in \ell^1_{+}$.  
\section{Proximal Calculus and Efficient Implementation} \label{sec:prox-calc}

The transformation in Section~\ref{sec:bfgs-sr1-metric} enables us to compute the proximal mapping with respect to the metric in the form of (\ref{metric}) via the proximal calculus developed in \cite{becker2019quasi}, which we state here for completeness. 
\begin{theorem}
    Let $B = B_0+U_1U_1^\top -U_2U_2^\top \in \opspd{\sigma}(\R^n)$ with $\sigma>0$, $B_0\in \opspd{\sigma}(\R^n)$ and $U_i\in\R^{n\times r_i}$ with rank $r_i$ ($i=1,2$). Set $B_1=B_0+U_1 U_1^\top $. Then, the following holds:
    \begin{equation}
        \prox{}{g}^B(x) = \prox{}{g}^{B_0}(x+B_1^{-1}U_2\alpha_2^*-B_0^{-1}U_1\alpha_1^*)\,,
    \end{equation}
    where $\alpha_i^*,i=1,2,$ are the unique zeros of the coupled system $\mathcal{L}(\alpha)=\mathcal{L}(\alpha_1,\alpha_2)=0$, where $\alpha=(\alpha_1,\alpha_2)\in \R^{r_1+r_2}$ and $\mathcal{L}=(\mathcal{L}_1,\mathcal L_2)$ is defined by
    \begin{equation}\label{eq:prox-calc-root-finding-prob}
        \begin{split}
            \mathcal L_1(\alpha_1,\alpha_2)&=U_1^\top (x+B_1^{-1}U_2\alpha_2- \prox{g}{}^{B_0}(x+B_1^{-1}U_2\alpha_2-B_0^{-1}U_1\alpha_1))+\alpha_1\,,\\
            \mathcal L_2(\alpha_1,\alpha_2)&=U_2^\top (x-\prox{}{g}^{B_0}(x+B_1^{-1}U_2\alpha_2-B_0^{-1}U_1\alpha_1)) +\alpha_2\,.\\
        \end{split}
    \end{equation}
    Here, $\map{\mathcal{L}}{\R^{r_1+r_2}}{\R^{r_1+r_2}}$ is Lipschitz continuous.
\end{theorem}
The computation of a possibly complicated proximal mapping $\prox{}{g}^B(x)$ is reduced to a simple (by assumption) proximal mapping $\prox{}{g}^{B_0}$ and a low dimensional root finding problem in \eqref{eq:prox-calc-root-finding-prob}, which we tackle by the semi-smooth Newton solver proposed in \cite{becker2019quasi}, formulated in Algorithm~\ref{Alg:Semismooth}. It requires to solve Newton-like equations where the classic Jacobian at the current iterate $\alpha\piter\k$ is replaced by the Clarke Jacobian $\partial^c \mathcal{L}(\alpha\piter\k)$ (see \cite{Clarke90}), where we account for inexact solutions of these subproblems in terms of an error $e_k$. 
\begin{algorithm}[th]
\caption{Semi-smooth Newton method to solve $\mathcal{L}(\alpha)=0$ in \eqref{eq:prox-calc-root-finding-prob}}\label{Alg:Semismooth}
\begin{algorithmic}

\Require [initial data] $\alpha_0\in\R^r$, [maximum iterations] $N$
\State \textbf{Update for $k = 0,\cdots,N$}:
\begin{enumerate}
    \item[(i)]
    Select $G_k\in\partial^c \mathcal{L}(\alpha_k)$, compute $\alpha_{k+1}$ such that
    \[\mathcal{L}(\alpha_k)+G_k(\alpha_{k+1}-\alpha_k)=e_k\,,\]
    and $e_k\in\R^r$ is an error term satisfying $\norm[]{e_k}\leq \rho_k\norm[]{G_k}$ and $\rho_k\geq 0$.
    \item[(ii)] \textbf{if} {$\mathcal{L}(\alpha_k)=0$} \textbf{then} \textbf{terminate.} 
\end{enumerate}
\State \textbf{End of for-loop}
\end{algorithmic}
\end{algorithm}
For completeness, we also state the convergence result of \cite{becker2019quasi} for Algorithm~\ref{Alg:Semismooth}. 
\begin{theorem}
If $g$ is in addition a tame function, then the Lipschitz continuous function $\mathcal{L}$ is semi-smooth \cite{bolte2009tame} and all elements of $\partial^C\mathcal{L}(\alpha^*)$ are non-singular\cite{becker2019quasi}. Therefore, if $\rho_k\leq \bar \rho$, $\forall k\in\N$, for some sufficiently small $\bar \rho$ and $\alpha_0$ sufficiently close to $\alpha^*$, then the sequence generated by the Algorithm~\ref{Alg:Semismooth} is well-defined and converges to $\alpha^*$ linearly. Additionally, if $\rho_k\to0$, the convergence is superlinear.
\end{theorem}
The tameness assumption is extremely mild, as it includes basically any function that occurs in practical applications, by excluding pathological special cases. For example this class of functions comprises all semi-algebraic functions \cite{bolte2009tame}.

\section{Numerical Experiment} \label{sec:experiment}

\begin{figure}[th]
    \centering
    \tikz{
        \node (A)   {\includegraphics[width=0.53\linewidth]{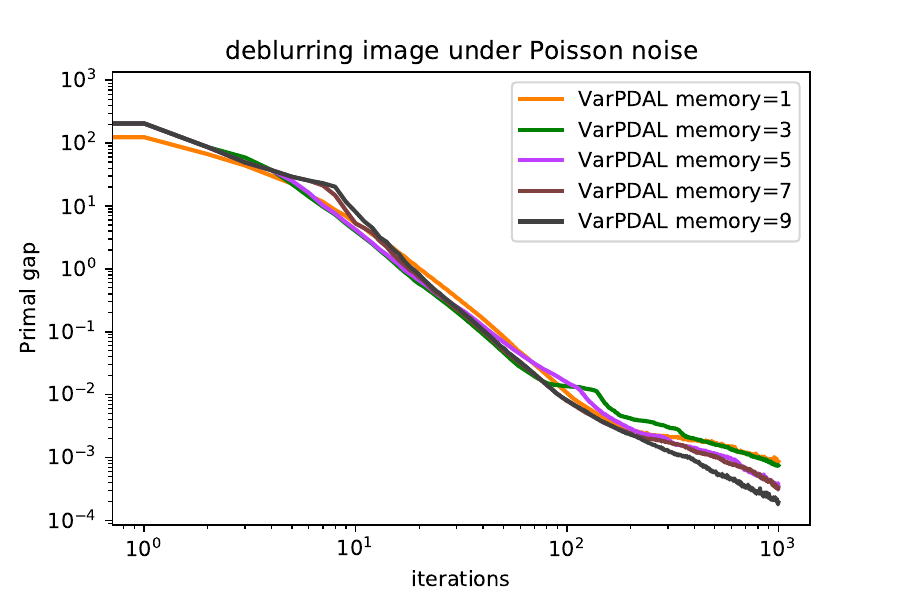}};
         \node[right=-1cm of A] (B) {\includegraphics[width=0.53\linewidth]{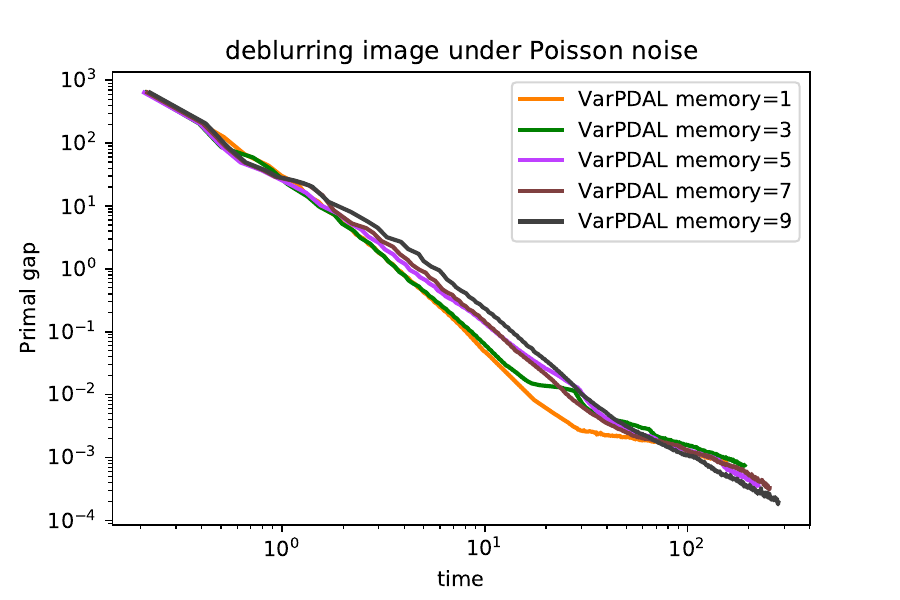}};
         \node[below= 1cm of A](C){
         \includegraphics[width=0.53\linewidth]{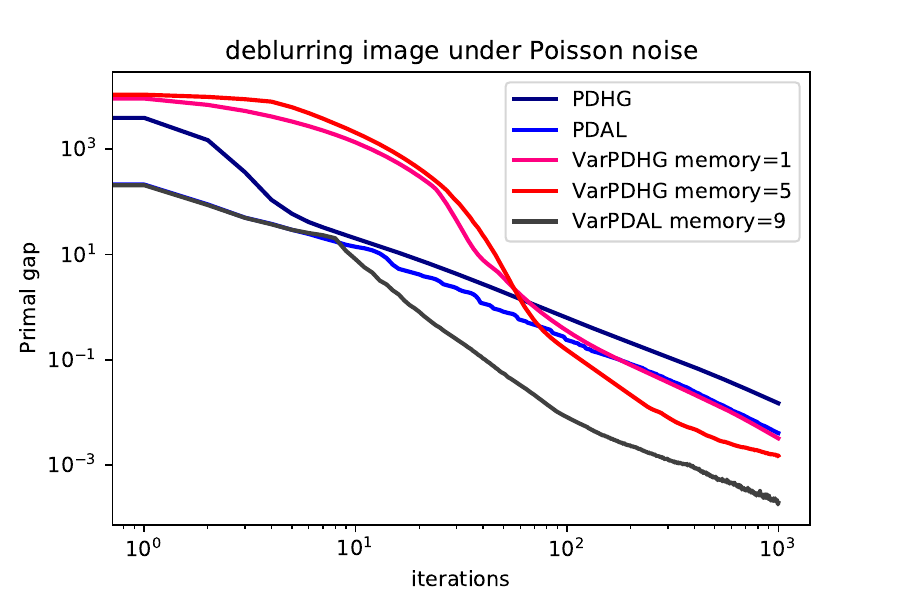}};
         
         \node[right = -1 cm of C](D){
         \includegraphics[width=0.53\linewidth]{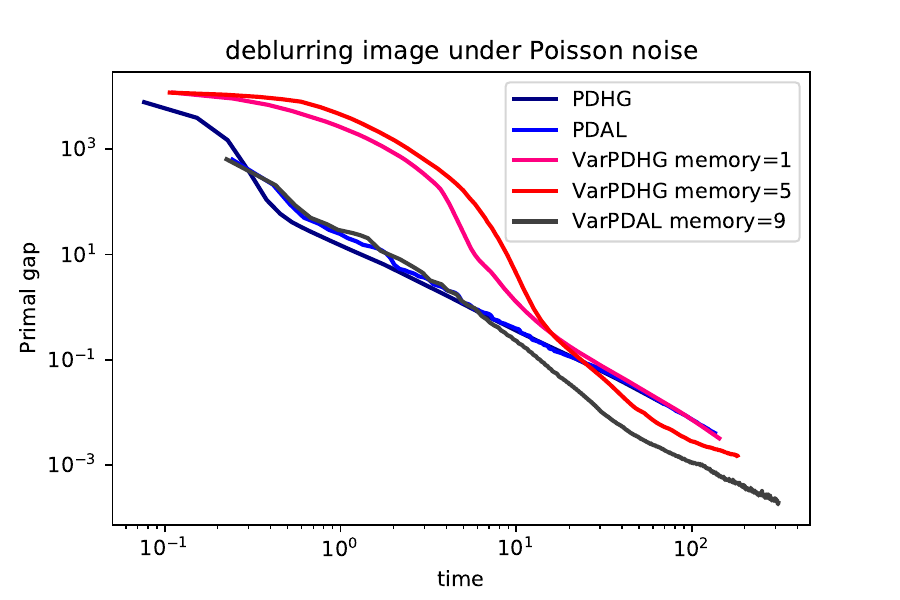
         }};
    }
    \caption{\label{fig:BFGSmemory}Performance evaluation for the experiment in \eqref{eq:num-exp-deblurr-poisson-saddle}. Algorithms are described in the text. Our algorithms, which combine a quasi-Newton variable metric with line search, outperform all other algorithms (that either use a variable metric \texttt{VarPDHG} or line search \texttt{PDAL}; or none of the two \texttt{PDHG}).}
\end{figure}

We apply our proposed algorithm for solving a challenging non-smooth image deblurring problem under a Poisson noise assumption \cite{vardi1985statistical}. Given the observation $b\in \R^{n_x\times n_y}$ as $n_x\times n_y$-sized image, the task is the following popular problem:
\begin{equation} \label{eq:deblurring-experiment-primal}
    \min_{x\in \mathbb{R}_+^{n_x \times n_y}} D_{KL}(b,Ax) + \gamma \norm[2,1]{\mathcal{D}x} \,,
\end{equation}
which involves the Kullback--Leibler divergence as data fidelity measure 
$
h(x)\coloneqq D_{\mathrm{KL}}(b,Ax)\coloneqq\sum_{i,j} (Ax)_{i,j}-b_{i,j} \log((Ax)_{i,j})
$ 
with respect to the blurred reconstruction $Ax$ with known blur operator and a discrete total variation regularization term $f(x)=\gamma \norm[2,1]{\mathcal{D}x}$ that is steered by a weight $\gamma>0$ where $\mathcal{D}$ implements discrete spatial finite differences. 
We recast \eqref{eq:deblurring-experiment-primal} into the saddle point problem:
\begin{equation}\label{eq:num-exp-deblurr-poisson-saddle}
    \min_{x\in \mathbb{R}^{n_x \times n_y}} \max_{y\in\R^{2\times n_x\times n_y}}
    \scal{\mathcal Dx}{y}+ \delta_{\mathbb{R}_+^{n_x \times n_y}}(x)+ D_{KL}(b,Ax) -\delta_{\norm[2,\infty]{\cdot}\leq \gamma} (y)
\end{equation}
and set $g(x):=\delta_{\mathbb{R}_+^{n_x \times n_y}}(x)$ and $K=\mathcal D$ in \eqref{eq:saddle-point-problem}. 
Figure~\ref{fig:BFGSmemory} compares several methods including PDHG with fixed stepsize (\texttt{PDHG}), PDHG with line search (\texttt{PDAL}), PDHG with fixed stepsize and variable metric (\texttt{VarPDHG}), PDHG with variable metric and line search (\texttt{VarPDAL}). The variable metric is generated by the limited memory BFGS method in Section~\ref{sec:bfgs-sr1-metric}. Figure~\ref{fig:BFGSmemory} shows the primal gap where the optimal primal value was computed for 10000 iterations by running \texttt{PDHG}. For the update of the variable metric \eqref{eq:metric-scalling}, we set $\gamma_1 =1.0$ and $\gamma_2 =0.99$. However, the Assumption~\ref{assumption1} is not satisfied since it is not guaranteed by the construction of the metric that there is a sequence $(\eta_k)_k\in\ell^1_+(\N)$ such that $M_{k+1}\preceq (1+\eta_k)M_k$. Fortunately, we still observe convergence of PDHG with this variable metric. Figure~\ref{fig:BFGSmemory} shows that a variable metric (\texttt{VarPDHG}) improves the convergence vs only using line search. However, our algorithm \texttt{VarPDAL} that combines both features is even faster, with the best performance when $m=9$. 

\begin{figure}[th]
    \centering
    \tikz{
        \node (A)   {\includegraphics[width=0.53\linewidth]{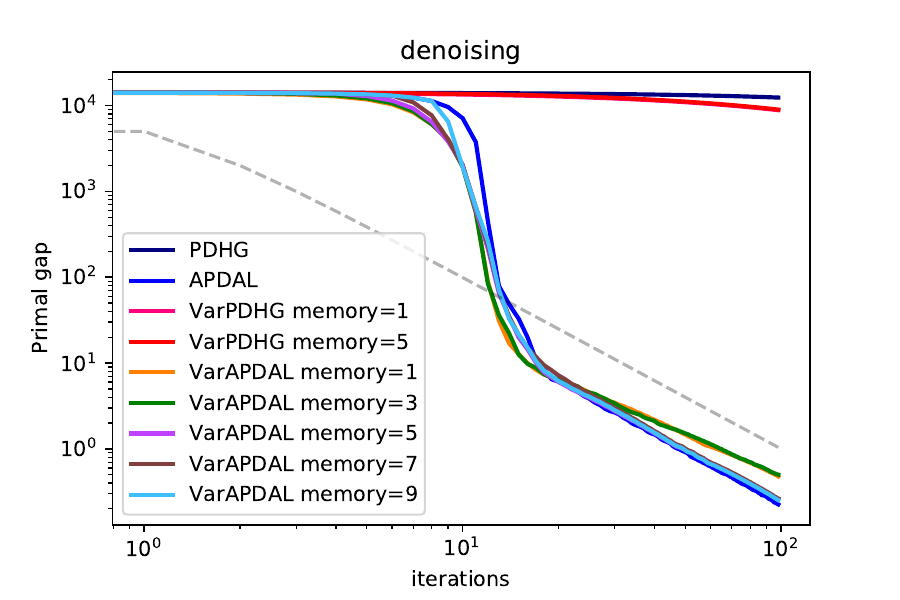}};
         \node[right=-1cm of A] (B) {\includegraphics[width=0.53\linewidth]{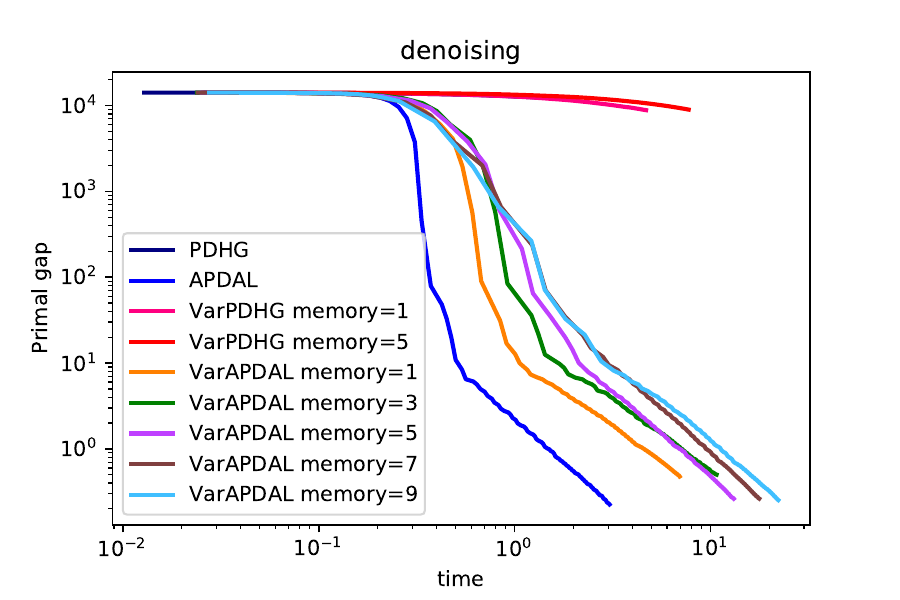}};
    
    }
    \caption{\label{fig:BFGSmemory-g-StrongConvex}Performance evaluation for the experiment for a strongly convex case. This figure reflects that our algorithms can retrieve $O(1/N^2)$ convergence rate as \texttt{PDAL}, which is theoretically guaranteed by Theorem \ref{thm:conv-rate-strong-convex}.}
\end{figure}

As a second experiment to test out accelerated version, we consider $A=\opid$ and a constraint set $\mathcal{C}\coloneqq \{x\vert x_{ij}\in [\epsilon,255]\}$ since the grey value of each pixel should be less than 255 and be positive; We set $\epsilon=0.1$.
In this case $D_{KL}(b,x)$ restricted on $\mathcal{C}$ is a strongly convex function. We apply the accelerated version of Algorithm \ref{Alg:QN-PDHG-line} in Theorem~\ref{thm:conv-rate-strong-convex}, and we use the notation \texttt{VarAPDAL}. Similary, \texttt{APDAL} denotes the accelerated version of \texttt{PDAL}. The dashed line in \ref{fig:BFGSmemory-g-StrongConvex} corresponds to $O(1/N^2)$. We can observe that \texttt{VarPDAL}, \texttt{APDAL} can converge faster than $O(1/N^2)$ as predicted by the convergence theorem.

In order to record the exact algorithmic details for our experiments, all code for experiments from this paper is available at \url{https://github.com/wsdxiaohao/VarPDAL.git}.

\section{Conclusion}

In this paper, we introduced a line search variant of a recently introduced quasi-Newton primal-dual algorithm. In contrast to related work, the employed quasi-Newton metric is of type ``identity $\pm$ low rank'', which captures significantly more second order information than a commonly used diagonal metric. We equally care for both, theoretical convergence guarantees including convergence rates as well as efficient practical implementation. The additional line search procedure usually leads to larger steps at a computational cost that pays off, which is confirmed by our numerical experiments. 



\hrule
\hrule
\hrule

\clearpage
\appendix
\section{Preliminaries}

There are several preliminaries we will use in the following section. The first one is a convergence result from \cite[Lemma 2.2.2]{Polyak} of a special sequence which appears in \ref{sec:proof-main-conv}.
\begin{lemma}\label{lemmaPolyak}
    Let $a_k\geq 0$ and let
    \begin{equation}
    \begin{split}
        &a_{k+1}\leq (1+\nu_{k})a_k + \zeta_k,\quad \nu_k\geq0,\ \zeta_k\geq 0, \\
        &\sum_{k\in\mathbb{N}}\nu_k<\infty,\quad\sum_{k\in\mathbb{N}}\zeta_k<\infty.\\
    \end{split}
    \end{equation}
    Then, $a_k\to A\geq 0$ for some $A<+\infty$.
\end{lemma}
The following identity is called (cosine rule), which proves to be very useful.
\begin{equation}\label{eq:quad}
    2\scal{a-b}{c-a} =\norm[]{b-c}^2 - \norm[]{a-b}^2 -\norm[]{a-c}^2\quad \forall a,b,c\in X\,.
\end{equation}
Another inequality appears many times in \ref{sec:proof-main-conv} is the characteristic property of the proximal operator with respect to a symmetric positive definite matrix $M$:
\begin{equation}\label{ineq:prox}
\hat x = \prox{g}{}^M(\bar x) \iff \scal{\hat x - \bar x}{y - \hat x}_{M} \geq g(\bar x) - g(y) \quad \forall y\in X\,.
\end{equation}
If $M=\opid$ is an identity matrix, then (\ref{ineq:prox}) is the characteristic property of the standard proximal operator.
Assume $(\hat x,\hat y)$ is a saddle point which solves (\ref{eq:saddle-point-problem}). Then we obtain
\begin{equation}\label{def:PrimalorDual-gap}
    \begin{split}
        P_{\hat x,\hat y}(x) &= g(x) + h(x) - g(\hat x)- h(\hat x) + \scal{K^*\hat y}{x-\hat x}\geq 0\quad \forall x\in X \,,\\
        D_{\hat x,\hat y}(y) &= f^*(y)-f^*(\hat y) - \scal{K\hat x}{y-\hat y}\geq 0 \quad    \forall y \in Y\,,\\
    \end{split}
\end{equation}
where $P_{\hat x,\hat y}(x)$ and $D_{\hat x,\hat y}(y)$ are convex.
Then $\mathcal G_{\hat x,\hat y}(x,y):= P_{\hat x,\hat y}(x)+D_{\hat x,\hat y}(y)$ is the primal-dual gap. Without ambiguity, in the proofs, we may omit the subscript in $P$ and $D$.

\section{Collection of Proofs}
\subsection{Proof of Lemma~\ref{lemma:sigma}}\label{sec:proof-lemma}
It is a similar argument with the one in \cite{malitsky2018first}.
\begin{enumerate}
    \item[(i)\&(ii)]$\sigma_k$ is decreased by $\mu\in(0,1)$ and the inequality (\ref{CB}) is satisfied as long as $\sigma_k <\underline{\sigma}_k\coloneqq\frac{-1+\sqrt{(4\delta\alpha)/\beta_k+1}}{  2\hat L} $ where $\hat L =\max \{L,L_K\}$.  We introduce a notation $\underline{\sigma}\coloneqq \frac{-1+\sqrt{(4\delta\alpha)/\beta+1}}{ 2\hat L} $. Since $\beta_k < \beta$, we have $\underline{\sigma}_k\geq\underline{\sigma}$. We argument by induction. We assume $\sigma_0 > \mu\underline{\sigma}_0$ and $\sigma_{k-1}>\mu\underline{\sigma}_{k-1}$ . For the case $\sigma_k=\bar \sigma_k$, then $\sigma_k\geq(\frac{\beta_{k-1}}{\beta_k})\sigma_{k-1}>\mu(\frac{\beta_{k-1}}{\beta_k})\underline{\sigma}_{k-1}>\mu\underline{\sigma}_k>\mu\underline{\sigma}$.  For the case $\sigma_k=\mu^i\bar \sigma_k$ ,  $\sigma_k'= \mu^{i-1}\bar \sigma_k$ does not satisfy (\ref{CB}). It follows $\sigma_{k}'>\underline{\sigma}_k $. Thus, $\sigma_k=\mu\sigma_k' >\mu \underline{\sigma}_k\geq\mu \underline{\sigma} $.
    \item[(iii)] By $\sigma_k\leq \sigma_{k-1}\sqrt{1+\theta_{k-1}}$, we get $\theta_k\leq \sqrt{1+\theta_{k-1}}$. Thus, $\theta_k$ is bounded from above.
    \hfill \qedsymbol
\end{enumerate}

\subsection{Proof of Theorem~\ref{thm:main-conv}} \label{sec:proof-main-conv}

The following proof is adapted from \cite{malitsky2018first}.
Assume $(\hat x,\hat y)$ is a saddle point of problem \ref{eq:saddle-point-problem} and $\beta_k\equiv \beta$. By using (\ref{ineq:prox}), we obtain the following two inequalities:
\begin{equation}\label{ineq:dual}
\scal{y^{k+1}-y^k -\sigma_k K x^{k+1}}{\hat y - y^{k+1}}\geq \sigma_k (f^*(y^{k+1}) -f^*(\hat y ))
\end{equation}
\begin{equation}
\scal{x^{k+1}- x^k +\tau_k M_k^{-1}K^*\bar y ^k+\tau_k M_k^{-1}\nabla h(x^k)}{\hat x - x^{k+1}}_{M_k} \geq \tau_k(g(x^{k+1})-g(\hat x))
\end{equation}
By using $\tau_k = \beta \sigma_k$
\begin{equation}\label{ineq:primal}
\begin{split}
    &\scal{\frac{1}{\beta}(x^{k+1}- x^k) + \sigma_k M_k^{-1}K^*\bar y ^k+\sigma_k M_k^{-1}\nabla h(x^k)}{\hat x - x^{k+1}}_{M_k} \\&\geq \sigma_k(g(x^{k+1})-g(\hat x))\\
\end{split}
\end{equation}
Similarly, we apply (\ref{ineq:prox}) on $y^k$ and obtain 
\begin{equation}
    \scal{y^k-y^{k-1} - \sigma_{k-1} K x^k}{y-y^k}\geq \sigma_{k-1}(f^*(y^k) -f^*(y)) \quad \forall y\in Y\,.
\end{equation}
Setting $y=y^{k+1}$ and $y = y^{k-1}$ respectively, we obtain
\begin{equation}\label{ineq:dual_k}
    \scal{y^k-y^{k-1} - \sigma_{k-1} K x^k}{y^{k+1}-y^k}\geq \sigma_{k-1}(f^*(y^k) -f^*(y^{k+1})) \quad \forall y\in Y\,,
\end{equation}
\begin{equation}\label{ineq:dual_km1}
    \scal{y^k-y^{k-1} - \sigma_{k-1} K x^k}{y^{k-1}-y^k}\geq \sigma_{k-1}(f^*(y^k) -f^*(y^{k-1})) \quad \forall y\in Y\,.
\end{equation}
We deduce from $\mathrm{(\ref{ineq:dual_k})}\times\theta_k$ and $\theta_k = \frac{\sigma_k}{\sigma_{k-1}}$ that:
\begin{equation}\label{ineq:theta1}
    \scal{\theta_k(y^k-y^{k-1})-\sigma_kKx^k}{y^{k+1}-y^k}\geq \sigma_k (f^*(y^k) - f^*(y^{k+1}))\,.
\end{equation}
By $ \mathrm{(\ref{ineq:dual_km1})}\times \theta_k^2$, we also get:
\begin{equation}\label{ineq:theta2}
    \scal{\theta_k(y^k-y^{k-1})-\sigma_k Kx^k}{\theta_k(y^{k-1}-y^k)}\geq \sigma_k(\theta_k f^*(y^k) -\theta_k f^*(y^{k-1}))\,.
\end{equation}
Summing (\ref{ineq:theta1}) and (\ref{ineq:theta2}) together, by using $\bar y ^k = y^k + \theta_k (y^k - y^{k-1})$, we obtain 
\begin{equation}\label{ineq:bar_y}
    \scal{\bar y^k - y^k-\sigma_k K x^k}{y^{k+1}-\bar y^k}\geq \sigma_k ((1+\theta_k)f^*(y^k)-\theta_kf^*(y^{k-1}) - f^*(y^{k+1}))\,.
\end{equation}
To sum up inequalties (\ref{ineq:dual}), (\ref{ineq:primal}) and (\ref{ineq:bar_y})
, we obtain 
\begin{equation}
\begin{split}
    &\scal{y^{k+1}-y^k -\sigma_k K x^{k+1}}{\hat y - y^{k+1}}\\&+ \scal{\frac{1}{\beta}(x^{k+1}- x^k) +\sigma_k M_k^{-1}K^*\bar y ^k+\sigma_k M_k^{-1}\nabla h(x^k)}{\hat x - x^{k+1}}_{M_k} \\&+ \scal{\bar y^k - y^k-\sigma_k K x^k}{y^{k+1}-\bar y^k}\\&\geq \sigma_k (f^*(y^{k+1}) -f^*(\hat y ))+\sigma_k(g(x^{k+1})-g(\hat x))+\sigma_k ((1+\theta_k)f^*(y^k)-\theta_kf^*(y^{k-1}) \\&- f^*(y^{k+1}))\,,
\end{split}
\end{equation}
Reorganizing the above inequality and using $\tau_k = \beta\sigma_k$, we have 
\begin{equation}\label{ineq:mix}
\begin{split}
    &\scal{y^{k+1}-y^k}{\hat y - y^{k+1}}+ \frac{1}{\beta}\scal{x^{k+1}- x^k }{\hat x - x^{k+1}}_{M_k} +\scal{\bar y^k - y^k}{y^{k+1}-\bar y^k} \\&+\scal{-\sigma_k K x^k}{y^{k+1}-\bar y^k}+ \scal{-\sigma_k K x^{k+1}}{\hat y - y^{k+1}}\\&+\scal{\sigma_k K^*\bar y ^k+\sigma_k\nabla h(x^k)}{\hat x - x^{k+1}} \\&\geq \sigma_k(g(x^{k+1})-g(\hat x))+\sigma_k ((1+\theta_k)f^*(y^k)-\theta_kf^*(y^{k-1}) - f^*(\hat y))\,,
\end{split}
\end{equation}
As in \cite{malitsky2018first}, we still have: 
\begin{equation}\label{eq:intermediate_term}
\begin{split}
    &\scal{-\sigma_k K x^k}{y^{k+1}-\bar y^k}+ \scal{-\sigma_k K x^{k+1}}{\hat y - y^{k+1}}+\scal{\sigma_k K^*\bar y ^k}{\hat x - x^{k+1}}\\
    &=\sigma_k\scal{Kx^k-Kx^{k+1}}{\bar y^k - y^{k+1}}+\sigma_k\scal{K\hat x}{\bar y^k -\hat y}-\sigma_k\scal{K^*\hat y}{x^{k+1}-\hat x}
\end{split}
\end{equation}
Adding $\sigma_k h(x^{k+1}) -\sigma_k h(\hat x)$ on both sides of (\ref{ineq:mix}), we obtain:
\begin{equation}\label{ineq:mix2}
\begin{split}
    &\scal{y^{k+1}-y^k}{\hat y - y^{k+1}}+ \frac{1}{\beta}\scal{x^{k+1}- x^k }{\hat x - x^{k+1}}_{M_k} +\scal{\bar y^k - y^k}{y^{k+1}-\bar y^k} \\&+\scal{-\sigma_k K x^k}{y^{k+1}-\bar y^k}+ \scal{-\sigma_k K x^{k+1}}{\hat y - y^{k+1}} \\&+\scal{\sigma_k K^*\bar y^k +\sigma_k\nabla h(x^k)}{\hat x - x^{k+1}} +\sigma_k h(x^{k+1}) -\sigma_k h(\hat x) \\&\geq \sigma_k(g(x^{k+1})-g(\hat x)+(1+\theta_k)f^*(y^k)-\theta_kf^*(y^{k-1}) - f^*(\hat y) +h(x^{k+1})-h(\hat x))\,.
\end{split}
\end{equation}
Combining (\ref{eq:intermediate_term}) and (\ref{ineq:mix2}), we have
\begin{equation}\label{ineq:mix3}
\begin{split}
    &\scal{y^{k+1}-y^k}{\hat y - y^{k+1}}+ \frac{1}{\beta}\scal{x^{k+1}- x^k }{\hat x - x^{k+1}}_{M_k} +\scal{\bar y^k - y^k}{y^{k+1}-\bar y^k} \\&\sigma_k\scal{Kx^k-Kx^{k+1}}{\bar y^k - y^{k+1}}+\sigma_k\scal{K\hat x}{\bar y^k -\hat y}-\sigma_k\scal{K^*\hat y}{x^{k+1}-\hat x}\\&+\scal{\sigma_k\nabla h(x^k)}{\hat x - x^{k+1}}+\sigma_k h(x^{k+1}) -\sigma_k h(\hat x) \\&\geq \sigma_k(g(x^{k+1})-g(\hat x)+(1+\theta_k)f^*(y^k)-\theta_kf^*(y^{k-1}) - f^*(\hat y) +h(x^{k+1})-h(\hat x))\,.
\end{split}
\end{equation}
By the definition of $D(y)$ (\ref{def:PrimalorDual-gap}) and $\bar y^k=y^k+\theta_k(y^k-y^{k-1})$, we have
\begin{equation}\label{eq:dualgap}
\begin{split}
    &(1+\theta_k)f^*(y^k)-\theta_kf^*(y^{k-1}) -f^*(\hat y)- \scal{K\hat x}{\bar y^k -\hat y} \\&= (1+\theta_k)(f^*(y^k)-f^*(\hat y )- \scal{K\hat x}{ y^k -\hat y})-\theta_k(f^*(y^{k-1})-f^*(\hat y)\\&-\scal{K\hat x}{ y^{k-1} -\hat y})\\&=(1+\theta_k)D(y^k) -\theta_k D(y^{k-1})\,.
\end{split}
\end{equation}
Using (\ref{eq:dualgap}) and the definition of $P(x)$, we deduce from (\ref{ineq:mix3}) that
\begin{equation}\label{ineq:mix4}
     \begin{split}
    &\scal{y^{k+1}-y^k}{\hat y - y^{k+1}}+ \frac{1}{\beta}\scal{x^{k+1}- x^k }{\hat x - x^{k+1}}_{M_k} +\scal{\bar y^k - y^k}{y^{k+1}-\bar y^k} \\&+\sigma_k\scal{Kx^k-Kx^{k+1}}{\bar y^k - y^{k+1}}+\scal{\sigma_k\nabla h(x^k)}{\hat x - x^{k+1}}\\&+\sigma_k h(x^{k+1}) -\sigma_k h(\hat x) \\&\geq \sigma_k(P(x^{k+1})+(1+\theta_k)D(y^k)-\theta_kD(y^{k-1}) )\,.
\end{split}
\end{equation}
From the line search condition (\ref{CB}), we have
\begin{equation}\label{ineq:cond}
\begin{split}
    &\sigma_k (h(x^{k+1})-h(x^k)-\scal{\nabla h(x^k)}{x^{k+1}-x^k})\\&\leq \frac{\delta}{2\beta}\norm[M_k]{x^{k+1}-x^k}^2 -\frac{1}{2}\sigma^2_k\norm[]{Kx^{k+1}-Kx^k}^2\,.
\end{split}
\end{equation}
Additionally, by the convexity of $h(x)$, we also have
\begin{equation}\label{ineq:convex h}
    h(x^{k})-h(\hat x) + \scal{\nabla h(x^k)}{\hat x -x^k}\leq 0\,. 
\end{equation}
Combining (\ref{ineq:cond}) and $\sigma_k\times\mathrm{(\ref{ineq:convex h})}$, we get
\begin{equation}\label{ineq:convex h2}
\begin{split}
    &\sigma_k (h(x^{k+1})-h(\hat x )-\scal{\nabla h(x^k)}{x^{k+1}-\hat x})\\&\leq \frac{\delta}{2\beta}\norm[M_k]{x^{k+1}-x^k}^2 -\frac{1}{2}\sigma^2_k\norm[]{Kx^{k+1}-Kx^k}^2\,.
\end{split}
\end{equation}
Thus, it follows from (\ref{ineq:mix4}) and (\ref{ineq:convex h2}) that
\begin{equation}
 \begin{split}
    &\scal{y^{k+1}-y^k}{\hat y - y^{k+1}}+ \frac{1}{\beta}\scal{x^{k+1}- x^k }{\hat x - x^{k+1}}_{M_k} +\scal{\bar y^k - y^k}{y^{k+1}-\bar y^k} \\&+\sigma_k\scal{Kx^k-Kx^{k+1}}{\bar y^k - y^{k+1}} + \frac{\delta}{2\beta}\norm[M_k]{x^{k+1}-x^k}^2 -\frac{1}{2}\sigma^2_k\norm[]{Kx^{k+1}-Kx^k}^2 \\&\geq \sigma_k(P(x^{k+1})+(1+\theta_k)D(y^k)-\theta_kD(y^{k-1}) )\,.
\end{split}
\end{equation}
Using Cauchy-Schwarz inequality, we obtain
\begin{equation}\label{ineq:mix5}
    \begin{split}
    &\scal{y^{k+1}-y^k}{\hat y - y^{k+1}}+ \frac{1}{\beta}\scal{x^{k+1}- x^k }{\hat x - x^{k+1}}_{M_k} +\scal{\bar y^k - y^k}{y^{k+1}-\bar y^k} \\&+ \frac{1}{2}\sigma^2_k\norm[]{Kx^k-Kx^{k+1}}^2+ \frac{1}{2}\norm[]{\bar y^k - y^{k+1}}^2 + \frac{\delta}{2\beta}\norm[M_k]{x^{k+1}-x^k}^2 \\&-\frac{1}{2}\sigma^2_k\norm[]{Kx^{k+1}-Kx^k}^2 \\&\geq \sigma_k(P(x^{k+1})+(1+\theta_k)D(y^k)-\theta_kD(y^{k-1}) )\,.
\end{split}
\end{equation}
Applying (\ref{eq:quad}), we deduce from (\ref{ineq:mix5}) that
\begin{equation}\label{ineq:conv-descent}
     \begin{split}
    &(\frac{1}{2}\norm[]{y^k-\hat y}^2 - \frac{1}{2}\norm[]{y^{k+1}-y^k}^2-\frac{1}{2}\norm[]{\hat y-y^{k+1}}^2) \\&+ (\frac{1}{2\beta}\norm[M_k]{x^{k}- \hat x }^2 - \frac{1}{2\beta}\norm[M_k]{x^{k+1}- x^k }^2-\frac{1}{2\beta}\norm[M_k]{ \hat x - x^{k+1}}^2)
     \\&+(\frac{1}{2}\norm[]{y^k- y^{k+1}}^2 - \frac{1}{2}\norm[]{\bar y^k - y^k}^2-\frac{1}{2}\norm[]{ y^{k+1}-\bar y^k }^2)\\& + \frac{1}{2}\norm[]{\bar y^k - y^{k+1}}^2 + \frac{\delta}{2\beta}\norm[M_k]{x^{k+1}-x^k}^2\\&\geq \sigma_k(P(x^{k+1})+(1+\theta_k)D(y^k)-\theta_kD(y^{k-1}) )\,.
\end{split}
\end{equation}
Reorganizing the above inequalities, we obtain 
\begin{equation}\label{ineq:conv-rate-relate}
    \begin{split}
    &\frac{1}{2}\norm[]{y^k-\hat y}^2 + \frac{1}{2\beta}\norm[M_k]{x^{k}- \hat x }^2 - \frac{1-\delta}{2\beta}\norm[M_k]{x^{k+1}- x^k }^2 \\&+\sigma_k \theta_kD(y^{k-1}) - \frac{1}{2}\norm[]{\bar y^k - y^k}^2
     \\&\geq \sigma_k(P(x^{k+1})+(1+\theta_k)D(y^k) ) +\frac{1}{2}\norm[]{\hat y-y^{k+1}}^2 +\frac{1}{2\beta}\norm[M_k]{ \hat x - x^{k+1}}^2\,.
\end{split}
\end{equation}
It follows from $\bar \sigma_k \leq \sqrt{1+\theta_{k-1}}\sigma_{k-1}$ that $\sigma_k \theta_k\leq \frac{\sigma_k^2}{\sigma_{k-1}}\leq \frac{\bar \sigma_k^2}{\sigma_{k-1}}\leq (1+\theta_{k-1})\sigma_{k-1}$.
Thus,
\begin{equation}
    \begin{split}
    &\frac{1}{2}\norm[]{y^k-\hat y}^2 + \frac{1}{2\beta}\norm[M_k]{x^{k}- \hat x }^2 - \frac{1-\delta}{2\beta}\norm[M_k]{x^{k+1}- x^k }^2 \\&+\sigma_{k-1} (1+\theta_{k-1})D(y^{k-1}) - \frac{1}{2}\norm[]{\bar y^k - y^k}^2
     \\&\geq \sigma_k(1+\theta_k)D(y^k)  +\frac{1}{2}\norm[]{\hat y-y^{k+1}}^2 +\frac{1}{2\beta}\norm[M_k]{ \hat x - x^{k+1}}^2\,.
\end{split}
\end{equation}
Since $(1+\eta_k)M_{k}\succeq M_{k+1}$, we can obtain the following key inequality:
\begin{equation}\label{ineq:descent}
    \begin{split}
    &\frac{1}{2}\norm[]{y^k-\hat y}^2 + \frac{1}{2\beta}\norm[M_k]{x^{k}- \hat x }^2 - \frac{1-\delta}{2\beta}\norm[M_k]{x^{k+1}- x^k }^2 \\&+\sigma_{k-1} (1+\theta_{k-1})D(y^{k-1}) - \frac{1}{2}\norm[]{\bar y^k - y^k}^2
     \\&\geq \sigma_k(1+\theta_k)D(y^k)  +\frac{1}{2}\norm[]{\hat y-y^{k+1}}^2 +\frac{1}{2\beta(1+\eta_k)}\norm[M_{k+1}]{ \hat x - x^{k+1}}^2\,.
\end{split}
\end{equation}
Set $A_k\coloneqq \frac{1}{2}\norm[]{y^k-\hat y}^2  + \sigma_{k-1}(1+\theta_{k-1})D(y^{k-1})+ \frac{1}{2\beta}\norm[M_{k}]{x^{k}-\hat x}^2$.
Then, we deduce from (\ref{ineq:descent}) that
\begin{equation}
    A_{k+1}\leq (1+\eta_k)A_k \,.
\end{equation}
By Lemma \ref{lemmaPolyak}, $ A_k$ is bounded from above by some constant $C$. Thus, $\norm[]{y^k-\hat y} $ and $ \norm[M_{k}]{x^{k}-\hat x}$ are both bounded. By the assumption that $M_k$ is uniformly bounded, $\norm[]{x^{k}-\hat x}$ is also bounded. As a result, we deduce from (\ref{ineq:descent}) that
\begin{equation}
\begin{split}
     \sum_{k}  \big(\frac{1-\delta}{2\beta}\norm[M_k]{x^{k+1}- x^k }^2 + \frac{1}{2}\norm[]{\bar y^k - y^k}^2 \big)& \leq \sum_{k} \big((1+\eta_k)A_k - A_{k+1}\big)\\&\leq C\sum_{k}\eta_k + A_0<+\infty\,.\\
\end{split}
\end{equation}
It implies that $\norm[M_k]{x^{k+1}- x^k }\to 0$ and $\norm[]{\bar y^k - y^k}\to 0$. So does $\norm[]{x^{k+1}- x^k }\to 0$, since $(M_k)_{k\in\N}\subset\mathcal{S}_\alpha(X)$. Since $\sigma_k>\sigma$ for some $\sigma$ which is shown in Lemma \ref{lemma:sigma} and $\beta>0$ is fixed, 
\begin{equation}
\begin{split}
    &\frac{y^{k+1}-y^{k}}{\sigma_{k}}=\frac{\bar y^{k+1}-y^{k+1}}{\sigma_{k+1}}\to 0 \quad \mathrm{as}\, k\to +\infty\,,\\&\frac{\norm[M_k]{x^{k+1}- x^k }^2}{\tau_k}\to 0\quad \mathrm{as}\, k\to +\infty\,.\\
\end{split}
\end{equation}
Since $(x^k,y^k)_{k\in \N}$ is bounded, we can extract a subsequence $(x^{k_i},y^{k_i})_{i\in \N}$ converging to some cluster point $(x^*,y^*)$. As in \cite{malitsky2018first}, similarly, by using the lower semi-continuity of functions $g$ and $f^*$ and the continuity of function $h$, we can pass the following two inequalities to the limit:
\begin{equation}
    \begin{split}
        &\scal{\frac{y^{k_i+1}-y^{k_i}}{\sigma_{k_i}} - K x^{k_i+1}}{ y - y^{k_i+1}}\geq  (f^*(y^{k_i+1}) -f^*( y ))\quad \forall y\in Y\,,\\
        &\scal{\frac{x^{k_i+1}- x^{k_i}}{\tau_{k_i}} + M_{k_i}^{-1}K^*\bar y ^{k_i}+M_{k_i}^{-1}\nabla h(x^{k_i})}{ x - x^{k_i+1}}_{M_{k_i}} \\&= \scal{\frac{M_{k_i}(x^{k_i+1}- x^{k_i})}{\tau_{k_i}}}{x-x^{k_i+1}}+\scal{ K^*\bar y ^{k_i}+\nabla h(x^{k_i})}{ x - x^{k_i+1}}_{}  \\&\geq (g(x^{k_i+1})-g( x))\quad \forall x\in X\,.\\
    \end{split}
\end{equation}
Thus, $(x^*,y^*)$ is the saddle point of \eqref{eq:saddle-point-problem}.
If, additionally, $f^*(y)\vert_{dom_{f^*}}$ is continuous, then $f^*(y^{k_i})\to f^*(y^*)$ and $D(y^{k_i})\to 0$ as $i\to +\infty$. 
From (\ref{ineq:descent}), we have 
$\frac{1}{\Pi^k_{j=1}(1+\eta_j)}A_k$ is monotone.
Setting $\hat x = x^*$ and $\hat y = y^*$ in (\ref{ineq:descent}), by the boundedness of $\sigma_k$ and $\theta_k$, it follows that
\begin{equation}
    \lim_{k\to \infty} \frac{A_k}{\Pi^\infty_{i=1}(1+\eta_i)}\leq\lim_{k\to \infty} \frac{A_k}{\Pi^k_{i=1}(1+\eta_i)} = \lim_{i\to \infty} \frac{A_{k_i}}{\Pi^{k_i}_{j=1}(1+\eta_j)}\leq \lim_{i\to \infty} A_{k_i}=0
\end{equation}
Since $\Pi^\infty_{i=1}(1+\eta_i)<+\infty$, we have $\lim_{k\to+\infty}A_k\to 0$ which means $x^k\to x^*$ and $y^k\to y^*$ as $k\to +\infty$. \hfill \qedsymbol

\subsection{Proof of Theorem~\ref{thm:conv-rate}} \label{sec:proof-conv-rate}

    We adapt the corresponding proof in \cite{malitsky2018first}. Let $\epsilon_k\coloneqq \sigma_k\big(P(x^{k+1})+(1+\theta_k)D(y^k)- \theta_k D(y^{k-1}) \big)$.
    Then we obtain the following inequality from (\ref{ineq:conv-descent}),
    \begin{equation}
        \frac{1}{2}\norm[]{y^k-\hat y}^2 -\frac{1}{2}\norm[]{y^{k+1}-\hat y }^2 +\frac{1}{2\beta}\norm[M_k]{x^k-\hat x}^2 -\frac{1}{2\beta}\norm[M_k]{x^{k+1}-\hat x}^2 - \frac{1}{2}\norm[]{\bar y^k - y^k}^2\geq \epsilon_k\,.
    \end{equation}
    By the assumption \ref{assumption1}, we get
    \begin{equation}\label{ineq:66}
        \frac{1}{2}\norm[]{y^k-\hat y}^2 -\frac{1}{2}\norm[]{y^{k+1}-\hat y }^2 +\frac{1}{2\beta}\norm[M_k]{x^k-\hat x}^2 -\frac{1}{2\beta}\frac{\norm[M_{k+1}]{x^{k+1}-\hat x}^2}{(1+\eta_k)} - \frac{1}{2}\norm[]{\bar y^k - y^k}^2\geq \epsilon_k\,.
    \end{equation}
    Since $(1+\eta_k)\geq 1$, it follows
    \begin{equation}
        \frac{1}{2}\norm[]{y^k-\hat y}^2 -\frac{1}{2}\frac{\norm[]{y^{k+1}-\hat y }^2}{(1+\eta_k)}+\frac{1}{2\beta}\norm[M_k]{x^k-\hat x}^2 -\frac{1}{2\beta}\frac{\norm[M_{k+1}]{x^{k+1}-\hat x}^2}{(1+\eta_k)}\geq \epsilon_k\,.
    \end{equation}
    Let both sides of the above inequality be divided by $\Pi_{i=1}^{k-1}(1+\eta_i)$ and it is common to assume that an empty product yields identity i.e. $\Pi_{i=1}^{0}(1+\eta_i)=1$.
    Thus,
    \begin{equation}\label{ineq:iterdecrease}
    \begin{split}
        &\frac{1}{2}\frac{\norm[]{y^k-\hat y}^2}{\Pi_{i=1}^{k-1}(1+\eta_i)}-\frac{1}{2}\frac{\norm[]{y^{k+1}-\hat y }^2}{\Pi_{i=1}^{k}(1+\eta_i)}+\frac{1}{2\beta}\frac{\norm[M_k]{x^k-\hat x}^2}{\Pi_{i=1}^{k-1}(1+\eta_i)} -\frac{1}{2\beta}\frac{\norm[M_{k+1}]{x^{k+1}-\hat x}^2}{\Pi_{i=1}^{k}(1+\eta_i)}\\&\geq \frac{\epsilon_k}{\Pi_{i=1}^{k}(1+\eta_i)}\,.
    \end{split}
    \end{equation}
    Summing up (\ref{ineq:iterdecrease}) for $k=1,\cdots,N$, we obtain 
    \begin{equation}
        \frac{1}{2}\norm[]{y^1-\hat y}^2+\frac{1}{2\beta}\norm[M_1]{x^1-\hat x}^2\geq \sum_{k=1}^{N}\frac{\epsilon_k}{\Pi_{i=1}^{k}(1+\eta_i)}\geq \sum_{k=1}^N \frac{\epsilon_k}{C}\,.
    \end{equation}
    Here, we used the $C=\sum_{k\in\N}(1+\eta_k)<+\infty$.
    
    The following steps are similar with the ones in \cite{malitsky2018first}. 
    \begin{equation}
    \begin{split}
         \sum_{k=1}^N\epsilon_k =& \sigma_N(1+\theta_N)D(y^k) +\sum_{k=2}^{N}[(1+\theta_{k-1})\sigma_{k-1}-\theta_k\sigma_k]D(y^{k-1})\\& - \theta_1\sigma_1D(y^0) + \sum_{k=1}^N\sigma_kP(x^{k+1})\,.\\
    \end{split}
    \end{equation}
    Since $D$ is convex, 
    \begin{equation}
        \begin{split}
            \sigma_N(1+\theta_N)D(y^N)+ &\sum_{k=2}^{N}[(1+\theta_{k-1})\sigma_{k-1}-\theta_k\sigma_k]D(y^{k-1})\\
            &\geq (\sigma_1\theta_1 + s_N)D(\frac{\sigma_1(1+\theta_1)y^1 + \sum_{k=2}^N\sigma_k\bar y^k}{\sigma_1\theta_1+s_N})\\
            &=(\sigma_k\theta_1 + s_N)D(\frac{\sigma_1\theta_1y^0+\sum_{k=1}^N\sigma_k \bar y^k}{\sigma_1\theta_1 +s_N})\\
            &\geq s_N D(\bar Y^N)\,, \\
        \end{split}
    \end{equation}
where $s_N = \sum_{k=1}^N \sigma_k$. Similarly,
\begin{equation}
    \sum_{k=1}^N \sigma_k P(x^{k+1}) \geq s_N P(\frac{\sum_{k=1}^N \sigma_k x^{k+1}}{s_N})=s_N P(\bar{X}^{N})\,.
\end{equation}
As a result,
\begin{equation}
    \mathcal{G}(\bar X^N,\bar Y^N) = P(\bar X^N) + D(\bar Y^N) \leq \frac{C}{s_N}\big( \frac{1}{2\beta}\norm[M_1]{x^1-\hat x}^2 + \frac{1}{2}\norm[]{y^1-\hat y}^2+\sigma_1\theta_1D(y^0)\big)\,.
\end{equation}
\hfill \qedsymbol
\subsection{Proof of Theorem \ref{thm:conv-rate-strong-convex} }\label{sec:proof-conv-rate-strong-convex}
The proof is also adapted from \cite{malitsky2018first}.
From the update formula of $\beta_k$, it follows that $\beta_k$ is decreasing.
First, we are going to prove that $\theta_k$ is bounded from above.
It is not difficult but tedious.
We know that if there exists a $C\in\R_+$ s.t $\theta_k\leq C\sqrt{1+\theta_{k-1}}$ then $\theta_k$ is bounded. From this, it is sufficient to prove that $\frac{\beta_{k-1}}{\beta_k}$ is uniformly bounded from above by some $C_\theta$.
According to 
 \begin{equation}
      \beta_k=\frac{\beta_{k-1}}{\min\{1+\frac{\gamma}{C_M}\beta_{k-1}\sigma_{k-1},C_\theta\}}\,,\ \forall \k\in\N\,, \quad \text{and}\quad \beta_0 >0 \,,
  \end{equation}
  we have that $\frac{\beta_{k-1}}{\beta_k}=\min\{1+\frac{\gamma}{C_M}\beta_{k-1}\sigma_{k-1},C_\theta\}\leq C_\theta$.

Second part, we are going to show the convergence rate.
Since $g$ is strongly convex, we obtain:
\begin{equation}\label{ineq:change}
\begin{split}
    &\scal{\frac{x^{k+1}- x^k}{\tau_k} + M_k^{-1}K^*\bar y ^k+ M_k^{-1}\nabla h(x^k)}{\hat x - x^{k+1}}_{M_k} \\&\geq (g(x^{k+1})-g(\hat x))+\frac{\gamma}{2}\norm[]{x^{k+1} -\hat x}^2\,.
\end{split}
\end{equation}
From Assumption \ref{assumption1}, it follows that for any $k\in\N$, 
\begin{equation}\label{eq:change metric}
    \frac{\gamma}{2}\norm[]{x^{k+1} -\hat x}^2\geq \frac{\gamma}{2C_M}\norm[M_{k+1}]{x^{k+1} -\hat x}^2\,.
\end{equation}

Following the same way in which we got equation \eqref{ineq:conv-rate-relate}, by equation \eqref{ineq:change} and the assumption that $ (1+\eta_k)M_k\succeq M_{k+1}$, we obtain 
\begin{equation}
\begin{split}
\frac{1}{2}\norm[]{y^k-\hat y}^2& -\frac{1}{2}\norm[]{y^{k+1}-\hat y }^2 +\frac{1}{2\beta_k}\norm[M_k]{x^k-\hat x}^2 - \frac{1-\delta}{2\beta_k}\norm[M_k]{x^{k+1}- x^k }^2 \\&-\frac{1}{2\beta_k}\frac{\norm[M_{k+1}]{x^{k+1}-\hat x}^2}{(1+\eta_k)} - \frac{1}{2}\norm[]{\bar y^k - y^k}^2\geq \epsilon_k+\frac{\gamma\sigma_k}{2}\norm[]{x^{k+1}-\hat x}^2\,.\\ 
\end{split}
\end{equation}

In order to obtain the following inequality, it is sufficient to assume $\delta\leq 1$. 
Thus,
\begin{equation}
\begin{split}
\frac{1}{2}\norm[]{y^k-\hat y}^2& -\frac{1}{2}\norm[]{y^{k+1}-\hat y }^2 +\frac{1}{2\beta_k}\norm[M_k]{x^k-\hat x}^2 \\&-\frac{1}{2\beta_k}\frac{\norm[M_{k+1}]{x^{k+1}-\hat x}^2}{(1+\eta_k)} - \frac{1}{2}\norm[]{\bar y^k - y^k}^2\geq \epsilon_k+\frac{\gamma\sigma_k}{2}\norm[]{x^{k+1}-\hat x}^2\,.\\ 
\end{split}
\end{equation}

Since $\delta\leq 1$, 
by dividing the above inequality with $\sigma_k$, we have
\begin{equation}
\begin{split}
\frac{1}{2\sigma_k}\norm[]{y^k-\hat y}^2& -\frac{1}{2\sigma_k}\norm[]{y^{k+1}-\hat y }^2 +\frac{1}{2\tau_k}\norm[M_k]{x^k-\hat x}^2 \\&-\frac{1}{2\tau_k}\frac{\norm[M_{k+1}]{x^{k+1}-\hat x}^2}{(1+\eta_k)} - \frac{1}{2\sigma_k}\norm[]{\bar y^k - y^k}^2\geq \frac{\epsilon_k}{\sigma_k}+\frac{\gamma}{2}\norm[]{x^{k+1}-\hat x}^2\,, \\ 
\end{split}
\end{equation}
where, we used $\tau_k = \beta_k\sigma_k$.
By using \eqref{eq:change metric}, from the above inequality, we obtain that
\begin{equation}
\begin{split}
     \frac{1}{2\sigma_k}\norm[]{y^k-\hat y}^2 &-\frac{1}{2\sigma_k}\norm[]{y^{k+1}-\hat y }^2 +\frac{1}{2\tau_k}\norm[M_k]{x^k-\hat x}^2 -\frac{1}{2\tau_k}\frac{\norm[M_{k+1}]{x^{k+1}-\hat x}^2}{(1+\eta_k)} \\&- \frac{1}{2\sigma_k}\norm[]{\bar y^k - y^k}^2\geq \frac{\epsilon_k}{\sigma_k}+\frac{\gamma}{2 C_M}\norm[M_{k+1}]{x^{k+1}-\hat x}^2\,.\\
\end{split}
\end{equation}
It follows from the above inequality that
\begin{equation}\label{ineq:decent-strongly-convex}
    \begin{split}
     \frac{1}{2\sigma_k}\norm[]{y^k-\hat y}^2 -\frac{1}{2\sigma_k}\norm[]{y^{k+1}-\hat y }^2 &+\frac{1}{2\tau_k}\norm[M_k]{x^k-\hat x}^2  - \frac{1}{2\sigma_k}\norm[]{\bar y^k - y^k}^2\\&\geq \frac{\epsilon_k}{\sigma_k}+\frac{1+(1+\eta_k)\tau_k\gamma/C_M}{2\tau_k(1+\eta_k)}\norm[M_{k+1}]{x^{k+1}-\hat x}^2\,,\\
\frac{1}{2\sigma_k}\norm[]{y^k-\hat y}^2 -\frac{1}{2\sigma_k}\norm[]{y^{k+1}-\hat y }^2& +\frac{1}{2\tau_k}\norm[M_k]{x^k-\hat x}^2  - \frac{1}{2\sigma_k}\norm[]{\bar y^k - y^k}^2\\&\geq \frac{\epsilon_k}{\sigma_k}+\frac{\tau_{k+1}(1+\tau_k\gamma/C_M)}{\tau_k}\frac{\norm[M_{k+1}]{x^{k+1}-\hat x}^2}{2\tau_{k+1}(1+\eta_k)}\,,\\
\end{split}
\end{equation}
For convenience, we set $\tilde \gamma = \gamma/C_M$.
From the update step of $\beta_k$, it follows that

\begin{equation}
    \frac{\tau_{k+1}(1+\tilde\gamma\tau_k)}{\tau_k}\geq \frac{\tau_{k+1}\min\{C_\theta,(1+\tilde\gamma\tau_k)\}}{\tau_k}= \frac{\sigma_{k+1}}{\sigma_k}
\end{equation}

Set $B_k\coloneqq\frac{1}{2\tau_k}\norm[M_k]{x^k-\hat x}^2+\frac{1}{2\sigma_k}\norm[]{y^k-\hat y}^2$ and $\tilde B_k \coloneqq \frac{B_k}{\Pi_{i=1}^{k-1}(1+\eta_i)}$. 
From \eqref{ineq:decent-strongly-convex}, we have:
\begin{equation}
    \frac{\sigma_{k+1}}{\sigma_k(1+\eta_k)}B_{k+1} + \frac{\epsilon_k}{\sigma_k}\leq B_k - \frac{1}{2\sigma_k} \norm[]{\bar y^k - y^k}^2
\end{equation}
By dividing the above inequality by $\Pi_{i=1}^{k-1}(1+\eta_i) \geq 1$, we obtain
\begin{equation}
    \frac{\sigma_{k+1}}{\sigma_k}\tilde{B}_{k+1} + \frac{\epsilon_k}{\sigma_k\Pi_{i=1}^{k-1}(1+\eta_i)}\leq \tilde{B}_k - \frac{1}{2\sigma_k\Pi_{i=1}^{k-1}(1+\eta_i)} \norm[]{\bar y^k - y^k}^2
\end{equation}
By multiplying $\sigma_k$ on both sides, we have
\begin{equation}
    \sigma_{k+1}\tilde{B}_{k+1} + \frac{\epsilon_k}{\Pi_{i=1}^{k-1}(1+\eta_i)}\leq \sigma_k\tilde{B}_k - \frac{1}{2\Pi_{i=1}^{k-1}(1+\eta_i)} \norm[]{\bar y^k - y^k}^2\,.
\end{equation}
By Assumption \ref{assumption1}, $C=\Pi_{i\in\N}(1+\eta_i)<+\infty$,
we have
\begin{equation}\label{ineq:B_k}
    \sigma_{k+1}\tilde{B}_{k+1} + \frac{\epsilon_k}{C}\leq \sigma_k\tilde{B}_k - \frac{1}{2C} \norm[]{\bar y^k - y^k}^2\,.
\end{equation}
Summing up \eqref{ineq:B_k} from $k=1,\cdots,N$, we obtain
\begin{equation}
    \sigma_{N+1}\tilde{B}_{N+1} + \sum
_{k=1}^N\frac{\epsilon_k}{C}\leq \sigma_1\tilde{B}_1 - \frac{1}{2C} \sum_{k=1}^{N}\norm[]{\bar y^k - y^k}^2\,.
\end{equation}
Since $\sigma_k$ is bounded by some $\sigma$ for any $k\in\N$, 
$\Tilde{B}_k$ is bounded from above. Since $C=\Pi_{i\in\N}(1+\eta_i)<+\infty $, $B_k$ is also bounded from above. So, $y^k$ is also bounded with $\lim_{k\to\infty}\norm[]{\bar y^k-y^k}^2= 0$. Thus, using the similar argument and notations in the proof \ref{sec:proof-main-conv}, we retrieve the same key inequality as the one in \cite{malitsky2018first}:
\begin{equation}
\begin{split}
    \mathcal{G}(\bar X^N,\bar Y^N)&\leq \frac{C}{s_N}(\sigma_1B_1 + \theta_1\sigma_1P(x^0)),\\
    \norm[M_{N+1}]{x^{N+1}-\hat x}^2&\leq \frac{C\tau_{N+1}}{\sigma_{N+1}}(\sigma_1A_1+\theta_1\tau_1P(x^0))=C\beta_{N+1}\,,\\
\end{split}
\end{equation}
Using the same argument from \cite{malitsky2018first}, we know from \ref{sec:proof-lemma} that $\sigma_k$ is bounded by $\mu\underline{\sigma}_k=\mu(\frac{-1+\sqrt{(4\delta\alpha)/\beta_k+1}}{  2\hat L})$ where $\hat L= \max\{L,L_K\}$. We claim that there exists a constant $C_\beta$ such that, $\beta_{k}=C_\beta(1/k^2)$.
\begin{enumerate}
    \item[i] If $\alpha\delta/(\beta_k)\leq 1$, by $ \sigma_k\geq \mu \underline{\sigma}_k\geq \mu\underline{\sigma}$, we have
    
    \begin{equation}
        \beta_{k+1} = \frac{\beta_k}{\min\{C_\theta,1+\tilde \gamma \beta_k\sigma_k\}}\leq\frac{\beta_k}{\min\{C_\theta,1+\mu\underline{\sigma}\delta\alpha\tilde \gamma\}}\,.
    \end{equation}
    
    In this case, $\beta_k$ decreases linearly. Thus, $\beta_{k+1}\leq C_\beta/(k+1)^2$ for $k$ sufficiently large.
    \item[ii] If $\alpha\delta/(\beta_k)\geq 1$, then $ \sigma_k>\mu\underline{\sigma}_k>\frac{\mu}{2\hat L}\sqrt{\frac{\delta\alpha}{\beta_k}}$. Therefore, for $k$ large enough, we have
    
    \begin{equation}
        \beta_{k+1} = \frac{\beta_k}{\min\{C_\theta,1+\tilde \gamma \beta_k\sigma_k\}}\leq\frac{\beta_k}{\min\{C_\theta,1+\frac{\mu\sqrt{\delta\alpha}\tilde \gamma}{2\hat L}\sqrt{\beta_k}\}}=\frac{\beta_k}{1+\frac{\mu\sqrt{\delta\alpha}\tilde \gamma}{2\hat L}\sqrt{\beta_k}}\,.
    \end{equation}
    
    In this case, by induction $\beta_k\leq \frac{C_\beta}{k^2}$ for some constant $C_\beta>0$.
\end{enumerate}
From $\sigma_k>\mu\underline{\sigma}_k>\mu\underline{\sigma}$, we have
$s_N= \sum_{k=1}^N \sigma_k >\sum_{k=1}^N \underline{\sigma}_k >\sum_{k=1}^N O(k)\sim N^2$ since $\beta_k\leq C_\beta/k^2$ for $k$ sufficiently large.
Then, we conclude the results.
\hfill \qedsymbol


\bibliographystyle{plain} 
\bibliography{main}

\end{document}